%% file: parkingf.tex
\documentclass[12pt,a4paper]{article}

\usepackage[english,francais]{babel} 
\RequirePackage[latin1]{inputenc}
\RequirePackage[T1]{fontenc}
\RequirePackage{graphicx,hyperref,ifpdf,color}

\RequirePackage{amsmath,amssymb,amsfonts,bbm}

\newcommand{\sth}{{{{\mbox{${\mbox{\scriptsize{\bf$\theta$}}}$}}}}}

\newcommand{\bp}{{\bf p}}

\newcommand{\bs}{{\bf s}}

\newcommand{\bX}{{\bf X}}
\newcommand{\bC}{{\bf C}}
\newcommand{\bF}{{\bf F}}

\newcommand{\btheta}{{\mbox{\boldmath$\theta$}}}

\newcommand{\R}{\mathbb{R}}
\newcommand{\N}{\mathbb{N}}

\newcommand{\Z}{\mathbb{Z}}
\renewcommand{\P}{\mathbb{P}}

\newcommand{\E}{\mathbb{E}}

\newcommand{\D}{\mathbb{D}}
\newcommand{\T}{\mathbb{T}}

\newcommand{\TT} { {\cal T }}

\newcommand{\FF} { {\cal F }}

\def\build#1_#2^#3{\mathrel{
\mathop{\kern 0pt#1}\limits_{#2}^{#3}}}

\def\cq{$\hfill \square$}
\def\un{\underline}
\def\d{{\rm d}}

\def\eps{\varepsilon}

\def\ba{\begin{eqnarray*}}
\def\ea{\end{eqnarray*}}

\def\wt{\widetilde}

\def\proof{\noindent{\bf Proof. }}

\def\supp{{\rm supp\,}}

\newcommand{\ind}{\mathbbm{1}}

\newtheorem{thm}{Theorem}
\newtheorem{lmm}{Lemma}
\newtheorem{prp}{Proposition}

\usepackage{hyperref} 
\usepackage{verbatim,multicol}

\def\bertoin{{\href{http://www.proba.jussieu.fr/pageperso/bertoin.html}{Jean
Bertoin}}}
\def\miermont{{\href{http://mahery.math.u-psud.fr/~miermont/}{Gr\'egory Miermont}}}
\def\dma{{\href{http://www.dma.ens.fr/}{DMA}}}
\def\lpma{{\href{http://www.proba.jussieu.fr/}{LPMA}}}

\title{Asymptotics in Knuth's parking problem for caravans}
\author{\bertoin\thanks{\lpma, Universit\'e Paris 6, 175 rue du
Chevaleret, F-75013 Paris} \& \miermont\thanks{\dma, \'Ecole Normale
Sup\'erieure. 45, rue d'Ulm, 75230 Paris Cedex
05, France}}

\begin{document}

\selectlanguage{english}

\maketitle

\begin{abstract}
We consider a generalized version of Knuth's parking problem, in which caravans
consisting  of a random number of cars arrive at random on the unit circle.
Then each car turns clockwise until it finds a free space to park.  Extending a recent
work by Chassaing and Louchard \cite{chasslou99}, we
relate the asymptotics for the sizes of blocks formed by occupied spots with the dynamics of the
additive coalescent.  According to the  behavior of the caravan's size
tail distribution, several qualitatively different versions of eternal additive
coalescent are involved. 
\end{abstract}

\bigskip


\bigskip

\noindent{\bf M.S.C. code: }60F17, 60J25. 

\newpage

\section{Introduction}\label{sec:intro}

The original parking problem of Knuth can be stated as follows.  Consider a
parking lot with $n$ spaces, identified with the cyclic group
$\Z/n\Z$. Initially the parking lot is empty, and  $m\leq n$ 
cars in a queue arrive one by one. Car $i$ tries to park on a uniformly
distributed space $L_i$ among the $n$ possible, 
independently of other cars, but
if the space is already occupied, then it tries places labeled
$L_i+1,L_i+2,\ldots$ until it finally finds a free spot to park. As cars
arrive, blocks of consecutive occupied spots are forming. It appears that a 
{\em phase transition}  occurs at the stage where the  parking lot 
is almost full,
more precisely when the number of free spots is of order $\sqrt{n}$. Indeed, 
while the largest block of occupied spots is of order $\log m$ with high probability 
as long as $\sqrt{m}=o(n-m)$, a block of size approximately $n$ is present
(while the others are of order at most $\log n$) with high probability 
when $n-m=o(\sqrt{m})$. In the meanwhile, precisely when $n-m$ is of order 
$\lambda\sqrt{m}$ with $\lambda>0$, a clustering phenomenon occurs as
$\lambda$ decays. The behavior of this clustering process has been studied
precisely by Chassaing and Louchard  \cite{chasslou99}. 
It turns out that the
process of the relative sizes of occupied blocks is related to the so-called 
{\em standard additive coalescent} \cite{jpse96cmc,jpda98sac}.

The model originates from a problem in Computer Science: spaces in the 
parking lot should be
thought of as elementary memory spaces, each of which can be used to store
elementary data (cars). Roughly, our aim in this work is to investigate
the more general situation where one wants to store larger files, each
requiring several elementary memory spaces. In other words,
single cars are replaced by {\it caravans} of cars, i.e. several cars may
arrive simultaneously at the same spot.
In this direction, it will be convenient to consider a continuous
version of the problem, that goes as follows. Let $p_{1},  \ldots p_{m}$ be a sequence
of positive real numbers with sum $1$, and $s_{1}, \ldots, s_{m}$, $m$ distinct locations 
on the unit circle $\mathbb{T}:=\mathbb{R}/\mathbb{Z}$. 
Imagine that $m$ drops of paint with masses  $p_{1},\ldots ,p_{m}$, 
fall successively at locations $s_{1},\ldots, s_{m}$. Each time a drop of paint falls,
we brush it clockwise in such a way that the resulting painted portion of $\mathbb{T}$
is covered by a unit density of paint. So at each step the drop of paint is used to cover
a new portion of the circle and the total length of the painted part of the circle 
when $i\leq m$ drops have fallen is $p_{1}+\cdots+ p_{i}$.
In this setting, drops of paint play the role of caravans,
and the painted portion of the circle corresponds to occupied spots in the 
parking lot. 

\begin{figure}
\begin{center}
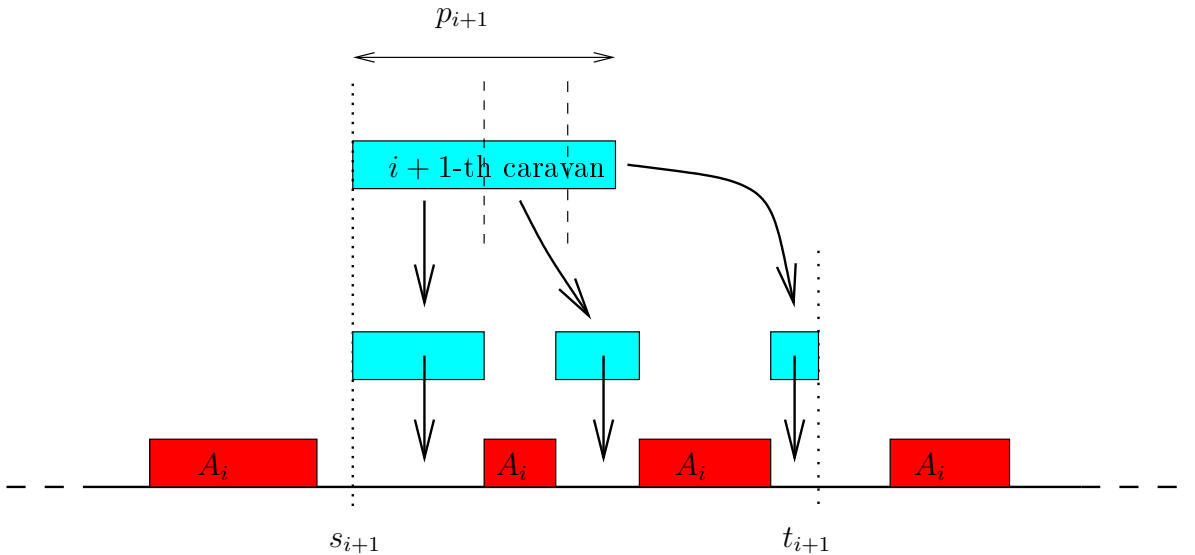
\end{center}
\caption{Arrival, splitting and parking of the $i+1$-th caravan in the process}
\label{fig:park}
\end{figure}

More precisely, we consider an increasing sequence $(A_0,\ldots,A_m)$ 
of open subsets of $\T$, starting from $A_0=\varnothing$ and ending at 
$A_m=\T$, 
which can be thought of as the successive painted portions of the circle. 
Given $A_i$ and the location $s_{i+1}$ from where the $i+1$-th drop of paint 
will be brushed, we paint as many space as possible to the right of $s_{i+1}$
with the quantity $p_{i+1}$ of paint, without covering the already painted 
parts, i.e.\ the blocks of $A_i$. 
Alternatively, we break the $i+1$-th caravan into several pieces, so 
that to fill as much as possible the holes left by $\mathbb{T}\setminus A_i$ 
after $s_{i+1}$, when 
reading in clockwise order. The last car to park 
arrives at some location $t_{i+1}$, and we let $A_{i+1}$ be the union of 
$A_i$ and the arc between $s_{i+1}$ and $t_{i+1}$, 
see Figure \ref{fig:park}. More formal definitions will come in 
Sect.\ \ref{sec:bridge}. 

In particular, $A_i$ is a disjoint union of open intervals and 
${\rm Leb}(A_i)=p_1+\ldots+p_i$. Let $\Lambda^{\bp}(i)
(=\Lambda(p_{1},\ldots,p_{m},
s_{1},\ldots,s_{m},i))$ be the sequence of the Lebesgue measures 
of the connected components of $A_{i}$, ranked by decreasing order.
It will be convenient to view $\Lambda^{\bp}(i)$ as an infinite sequence,
by completing with an infinite number of zero terms.

Now consider the following random problem. Let $\ell>0$ be a random variable 
with finite expectation $\mu_{1}=\E[\ell]$. We say that $\ell\in{\cal D}_{2}$ whenever
$\ell$ has a finite second moment $\mu_{2}=\E[\ell^2]$. For $\alpha\in(1,2)$, we say that 
$\ell\in{\cal D}_{\alpha}$ whenever 
\begin{equation}
\label{tail}
\P(\ell>x)\build\sim_{x\to\infty}^{} c x^{-\alpha}
\end{equation}
for some $0<c<\infty$. This implies that $\ell$ is in the domain of attraction of a spectrally
positive stable random variable with index $\alpha$, and we stress that our results 
can be extended
under this more general hypothesis; (\ref{tail}) is 
only intended to make things easier. We suppose from now on that
 $\ell\in{\cal D}_{\alpha}$ for some 
$\alpha\in(1,2]$, and take a random iid sample 
$\ell_{1},\ell_{2},\ldots$ of variables 
distributed as $\ell$, and independently of this sequence, 
 iid   uniform random variables on $[0,1)$, $U_{1},U_{2},\ldots$.
For $\eps>0$, set 
$$T_{\eps}=\inf\{i:\ell_{1}+\ldots+\ell_{i}\geq 1/\eps\},$$
so by the elementary renewal theorem, $T_{\eps}\sim 1/(\eps \mu_{1})$. Then
introduce the sequence $(\ell^*_{i},1\leq i\leq T_{\eps})$ defined by 
$$\ell^*_{i}=\ell_{i} \hbox{ for }1\leq i\leq T_{\eps}-1\ \hbox{ and } \ 
\ell^*_{T_{\eps}}=\eps^{-1}-(\ell_{1}+\ldots+\ell_{T_{\eps}-1}),$$ so the terms 
of $\ell^*$ sum to $1/\eps$. 

Following Chassaing and Louchard \cite{chasslou99}, we are interested in the formation of macroscopic painted components in the limit when
$\eps$ tends to $0$, at times close to $T_{\eps}$, i.e.  when the circle is almost entirely painted.
Specifically, we let 
$$\bX^{(\eps)}(t)=\Lambda^{\bp}(T_{\eps}-\lfloor 
t\eps^{-1/\alpha}\rfloor)\,,\qquad t\geq 0\,,$$ for  $\Lambda^{\bp}$
defined as
above with the data $m=T_{\eps},p_{i}={\eps}\ell^*_{i},s_{i}=U_{i}$. 
Observe that $T_{\varepsilon}-[t\varepsilon^{-1/\alpha}]$ decreases when
$t$ increases, and therefore, in order to investigate the formation of 
painted
components, we should consider the
process $\left({\bf X}^{(\varepsilon)}(t), t\geq0\right)$
backwards in time. This is what we shall do in Theorem 1, using the
exponential time change $t\to e^{-t}$.

Before describing our main result,  let us first recall some features of
the additive coalescent. The additive coalescent $\bC$ is a Markov process
with values in the infinite ordered simplex
$$S=\left\{\bs=(s_{1},s_{2},\ldots):s_{1}\geq s_{2}\geq \ldots\geq 0, 
\sum_{i=1}^{\infty} s_{i}\leq 1\right\}$$
endowed with the uniform distance,
whose evolution is described formally by: given that the current state is 
$\bs$, two terms $s_i$ and $s_j$, $i< j$, of $\bs$ are chosen and merge into
a single term $s_i+s_j$ (which implies some reordering of the resulting 
sequence)  at a rate  equal to  $s_i+s_j$. A  version $(\bC(t),t\in\R)$  of this
process defined for times describing the whole real axis is called {\em eternal}. 
We refer to \cite{jpda98sac, jpda97ebac} for background.
% The {\em standard} version is the unique (\cite{jpda97ebac}) one such that no 
% terms of $C(t)$ is preponderant in the scale $e^{-t}$ as $t\to-\infty$:
% \begin{equation}
% \label{asymC}
% e^{-t}C_i(t)\build\to_{t\to-\infty}^{}0 \qquad i\geq 1. 
% \end{equation}

As shown in \cite{bertfrag99}, eternal additive coalescents can be encoded by certain
bridges. Specifically, let $B=(B(x), 0\leq x \leq 1))$ 
be a c\`adl\`ag  real-valued process with exchangeable increments, such that
$B(0)=B(1)=0$. Suppose further that $B$ 
has infinite variation and  no negative jumps a.s. Then $B$ attains its overall infimum at a unique location 
$V$ (which is uniformly distributed on $[0,1]$), and $B$ is continuous at $V$.
Consider the so-called {\em Vervaat transform} which maps the bridge $B$ 
into an excursion ${\cal E}$ defined by
$${\cal E}(x)=B(V+x)-B(V) \quad \hbox{for  }  0\leq x\leq 1,$$
where the addition $V+x$ is modulo $1$. Finally, we let for $t\geq 0$
$${\cal E}^{(t)}_x=
{\cal E}(x)-tx\quad \hbox{for  }  0\leq x\leq 1,$$ and introduce
$\bF(t)$ as the random element of 
$S$ defined by the ranked sequence of the lengths of the constancy intervals of the process
$\un{{\cal E}}^{(t)}=(\inf_{0\leq y\leq x}{\cal E}^{(t)}(y),0\leq x\leq 1)$.  
Here, a constancy interval means a connected component of the complement of the support of the 
Stieltjes measure $\d(-\un{{\cal E}})$.  Finally, if we define $\bC(t)=\bF(e^{-t})$, then
$\bC=(\bC(t), -\infty < t < \infty)$ is an eternal additive coalescent (see Section \ref{subsec:extreme} 
for comments and details).

In this work, eternal additive coalescent associated to certain remarkable bridges will play 
a key role. More precisely, we write $\bC^{(2)}=(\bC^{(2)}(t), -\infty < t < \infty)$
for the eternal additive coalescent $\bC$ constructed above
when $B={\cal B}^{(2)}$ is a standard Brownian bridge; so that $\bC^{(2)}$ is the so-called 
standard additive coalescent (cf. \cite{bertfrag99, jpda98sac}). 
Next, for $1<\alpha < 2$, we denote by $\sigma^{(\alpha)}=(\sigma^{(\alpha)}(t), 
t\geq 0)$ a standard spectrally positive stable L\'evy process with index $\alpha$,
that is $\sigma^{(\alpha)}$ has independent and stationary increments, no negative jumps,
and
$$\E(\exp(-\lambda \sigma^{(\alpha)}(t)))=\exp(t\lambda^{\alpha})\,,
\qquad\hbox{for all }\lambda\geq 0.$$
We call {\em standard stable loop}
\footnote{We call ${\cal B}^{(\alpha)}$ a {\it loop} 
and not a {\it bridge} to avoid a possible confusion:
even though ${\cal B}^{(\alpha)}$ starts from $0$, ends at $0$ and
 has exchangeable increments, it does not have the same law as 
 the stable process $\sigma^{(\alpha)}$  conditioned on $\sigma^{(\alpha)}(1)=0$!} 
of index $\alpha$ the process
 ${\cal B}^{(\alpha)}$ defined by
\begin{equation}\label{pontstable}
{\cal B}^{(\alpha)}(x)=\sigma^{(\alpha)}(x)-x\sigma^{(\alpha)}(1),
\qquad \hbox{for }0\leq x \leq 1.
\end{equation}
We finally write  $\bC^{(\alpha)}=(\bC^{(\alpha)}(t), -\infty < t < \infty)$
for the eternal additive coalescent $\bC$ constructed above
when the bridge $B$ is the standard stable loop of index $\alpha$.

We are now able to state our main result.

\begin{thm}\label{convaddcoal}
The process $(\bX^{(\eps)}(t),0\leq  t< T_{\eps})$ converges as
$\eps\downarrow0$ in the sense of weak convergence of finite-dimensional 
distributions to some process $\bX=(\bX(t), 0\leq t < \infty)$.
The exponential time-changed process
$(\bX(e^{-t}), -\infty< t < \infty)$ is an eternal additive coalescent; more precisely:

\noindent {\rm (i)} When $\alpha=2$, $(\bX(e^{-t}), -\infty< t < \infty)$ is distributed as
$$(\bC^{(2)}(t+{1\over 2}\log (\mu_{2}/\mu_{1})-\log \mu_{1}), -\infty<t<\infty).$$

\noindent {\rm (ii)} When $1<\alpha<2$, $(\bX(e^{-t}), -\infty< t < \infty)$ 
is distributed as
$$(\bC^{(\alpha)}(t+{1\over \alpha}\log 
\left({\Gamma(2-\alpha)c\over (\alpha-1)\mu_{1}}\right)-\log \mu_{1}), -\infty<t<\infty).$$
\end{thm}

It might be interesting to discuss further the role of the parameter
$\alpha$ and the interpretation in terms of phase transition.
As it was already mentioned, the renewal theorem
entails than the number of drops of paint needed for the complete
covering is $T_{\varepsilon}\sim 1/(\varepsilon \mu_1)$, a quantity which
is not sensitive to $\alpha$. It is easy to show that for every $a<1$,
there are no macroscopic painted components when only $[aT_{\varepsilon}]$
drops of paint have fallen, so the phase transition (i.e. the number of
drops which is needed for the appearance of macroscopic components)
occurs for numbers close to
$T_{\varepsilon}$. More precisely, the regime for the phase transition is
of order $T_{\varepsilon} -\varepsilon^{-1/\alpha}$; so the phase
transition occurs closer to $T_{\varepsilon}$ when $\alpha$ is larger.
We would like also to stress that one-dimensional distributions of the
limiting additive coalescent process ${\bf X}$ depend on $\alpha$, but
not its semigroup which is the same for all $\alpha\in(1,2]$.
A heuristic explanation might be the following: the number of drops
needed to complete the covering once the phase transition has occurred
is too small (of order $\varepsilon^{-1/\alpha}$) to observe significant
differences in the dynamics of aggregation of macroscopic painted
components.
 
 \medskip
 
 \noindent{\bf Remark.} Our model bears some similarity with another 
parking problem on the circle,
where drops of paints fall uniformly on the circle and then are brushed 
clockwise, 
but where overlaps are now allowed (some points may be covered this way 
several times), 
call it the ``random covering of an interval'' problem. 
However, as showed in \cite{bertcov03}, this last model has
very different asymptotics from those of the parking problem, as it
turns out that the random covering of an interval is related to Kingman's
coalescent rather than the additive coalescent. A shared feature is that the
phase transition of the random covering problem appears also when the circle
is almost completely covered, but for example the different fragments are
ultimately finite in number rather than infinite.  

We also mention yet another parking problem, first considered by R\'enyi 
(see \cite{renpark,DRparking}). 
In can be formulated as follows: caravans with size 
$\eps$ are placed on $\T$ (the original work rather considers $(0,1)$) one 
after another, but the locations $s_i$ 
where cars park are chosen uniformly among
spaces that do not induce overlaps and splitting of caravans, i.e.\ so that 
the length of the arc from $s_i$ to $t_i$ is exactly $\eps$. This is done 
until no uncovered sub-arc of $\T$
with size $\geq \eps$ remains. This process does not involve coalescing 
blocks of cars, and one is rather interested in the properties of the 
random number of cars that are able to park. 

\medskip

  The method in \cite{chasslou99} relies on an encoding {\em parking} function 
which is shown to be asymptotically related to a function of standard Brownian
bridge, and a representation of the standard additive coalescent due to
Bertoin \cite{bertfrag99}.
Our approach to Theorem \ref{convaddcoal} is close in spirit to that of \cite{chasslou99},
and uses the representation of eternal additive 
coalescent that we presented above; we briefly sketch it here. 
First, we encode the process $\bX^{(\eps)}$ by a
bridge with exchangeable increments in Sect.\ \ref{sec:bridge}. In Sect.\
\ref{sec:conv},  we  show that  this  bridge  converges  to some  bridge  with
exchangeable increments  that can be represented  in terms of the standard 
Brownian bridge (for $\alpha=2$) or the standard stable loop (for $1<\alpha<2$). 
Theorem \ref{convaddcoal} then follows readily.

The rest of this work is organized as follows. In Section \ref{sec:bridge} we 
provide a representation
of the painted components in terms of a bridge and its Vervaat's transform. The convergence
of  these  bridges  when  $\eps$  tends  to  $0$  is  established  in  Section
\ref{sec:conv}, 
and that of the sequence
of the sizes of the painted components in Section \ref{convXeps}. Section \ref{sec:discrete}
is devoted to a brief discussion
of the analogous discrete setting (i.e. Knuth's parking for caravans), 
and finally some complements are presented in Section \ref{sec:compl}.

\section{Bridge representation}\label{sec:bridge}

We develop a representation of the parking process with the help of 
bridges with exchangeable 
increments, which is crucial to our study. 

Let us first give the proper definition the of sequence $(\varnothing=
A_0,\ldots,A_m)$ of the 
Introduction. We identify the circle $\mathbb{T}$ with $[0,1)$ and write
$p_{\mathbb{T}}: \mathbb{R}\rightarrow \mathbb{T}$ for
the canonical projection. 
If $A$ is a measurable subset of $\T$ (identified with $[0,1)$), let 
$F_A$ be its repartition function defined by $F_A(x)={\rm Leb}([0,x]\cap A)$ 
for $0\leq
x<1$, where ${\rm  Leb}$ is Lebesgue measure. Also, extend  $F_A$ on the whole
real line with the formula $F_A(x+1)=F_A(x)+F_A(1-)$. 
Given $A_i$ for some $0\leq i\leq m-1$, let 
$$t_{i+1}=\inf\{x\geq s_{i+1}:F_{A_i}(x)+p_{i+1}-(x-s_{i+1})
\leq F_{A_i}(s_{i+1})\}.$$ 
Notice that the arc $p_{\T}((s_{i+1},t_{i+1}))$ oriented clockwise from 
$s_{i+1}$ to $p_{\T}(t_{i+1})$ has length $t_{i+1}-s_{i+1}\geq p_{i+1}$. Then 
let $A_{i+1}$ be the interior of the closure of $p_{\T}((s_{i+1},t_{i+1}))\cup 
A_i$. 
%Thus, rather than
%simply adding the arc $p_{\T}((s_{i+1},s_{i+1}+p_{i+1}))$ to $A_i$, we break 
%it in several pieces $p_{T}((s_{i+1},t_{i+1}))\setminus  A_i$ and fit them in
%the ``holes'' left by $A_i$, starting at the
%point $s_{i+1}$. 
The point in taking the closure and then the interior is that we 
consider that two painted connected components of $\mathbb{T}$ 
that are at distance $0$ constitute in fact a single painted connected 
component. 

\begin{figure}
\begin{center}
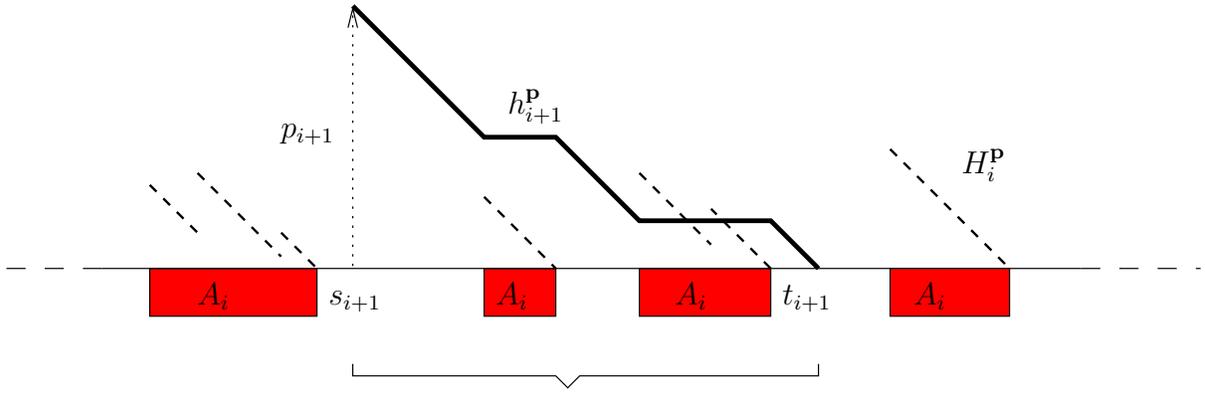
\end{center}
\caption{The function $h^{\bp}_{i+1}$ (thick line) corresponding to the 
$i+1$-th caravan of Figure \ref{fig:park}. The blocks of $A_i$ are represented
under the axis, and the dashed lines represent 
the profile $H^{\bp}_{i}$ (it gives more information than $A_i$ alone). The
bracket under the figure indicates how $A_{i+1}$ is obtained by 
formation of a new block comprising the blocks of $A_i$ between 
$s_{i+1}$ and $t_{i+1}$}
\label{fig:park2}
\end{figure}

Define
$$h^{\bp}_{i+1}(x)=F_{A_i}(x)-F_{A_i}(s_{i+1})+p_{i+1}-(x-s_{i+1}) 
\qquad s_{i+1}\leq x\leq t_{i+1},$$
and  $h^{\bp}_{i+1}(x)=0$ in $[t_{i+1},s_{i+1}+1)$, so $h^{\bp}_i$ is a
c\`adl\`ag function (right-continuous with left-limits) on $[s_{i+1},s_{i+1}+1)$.  Consider 
it as a function on $\T$ by letting $h^{\bp}_{i+1}(x)=h^{\bp}_{i+1}(y)$ where $y$ is the 
element of $[s_{i+1},s_{i+1}+1)\cap p_{\T}^{-1}(x)$. 
The quantity $h^{\bp}_{i+1}(x)$
can be thought of as the quantity of cars of the $i+1$-th caravan that try to park at $x$. See Figure \ref{fig:park2}.

We consider the {\em profile} 
\begin{equation}\label{def:profile}
H^{\bp}_i=\sum_{j=1}^i h^{\bp}_j
\end{equation}
of the parking at
step $0\leq i\leq m$, so $H^{\bp}_i(x)$ is the total quantity of cars that 
have tried (successfully or not) to park at $x$ (with the convention that
$H^{\bp}_i(1)=H^{\bp}_i(0)$) before the $i+1$-th caravan has arrived. 
\begin{lmm}\label{H}
For $1\leq i\leq m$,

\noindent
{\rm (i)} the set $A_i$ is the interior of the support of 
$H^{\bp}_{i}$. 

\noindent
{\rm (ii)} $H^{\bp}_i(t_i-)=0$. 

\noindent
{\rm (iii)} $H^{\bp}_i$ jumps at times $s_1,\ldots,s_i$
with respective jump magnitudes $p_1,\ldots,p_i$, and has a drift with slope
$-1$ on its support.That is, if $[v,v']\subseteq\supp(H^{\bp}_i)$, 
$$H^{\bp}_i(x+v)=H^{\bp}_i(v-)-x+\sum_{j=1}^ip_j\ind_{\{v\leq s_j\leq v+x\}}
\qquad 0\leq x\leq v'-v.$$
\end{lmm}
\proof
Properties (i) and (iii) are easily shown using a recursion on $i$ and splitting the
behavior of $h^{\bp}_i$ on $A_{i-1}$ and $A_i\setminus A_{i-1}$. We give some
details for (ii).  For $i\geq 1$, notice that by definition  $t_i$ cannot be a
point of increase of $F_{A_{i-1}}$, i.e. a point such that 
$F_{A_{i-1}}(t_i-\eps)<F_{A_{i-1}}(t_i)<F_{A_{i-1}}(t_i+\eps)$ for every
$\eps>0$.  Therefore, $t_i\notin A_{i-1}$ and
$h^{\bp}_j(t_i)=h^{\bp}_j(t_i-)=0$ for  $j<i$. Since it  follows by continuity
of $F_{A_{i-1}}$ that $h^{\bp}_i(t_i-)=0$, (ii) is proved. \cq

Consider the {\em bridge function}:
$$b^{\bp}_i(x)=-x+\sum_{j=1}^ip_j\ind_{\{x\geq s_j\}} \qquad 0\leq x<1,$$
which starts from $b^{\bp}_i(0)=0$ and ends at $b^{\bp}_i(1-)=p_{1}+\ldots + p_{i}-1$.
We extend $b^{\bp}_i$ to a function on $\R$ by setting $b^{\bp}_i(x+1)=b^{\bp}_i(x)+
b^{\bp}_i(1-)$. For any $v\in[0,1)$, it is easily seen using (iii) in Lemma
\ref{H} that 
$$H^{\bp}_i(x+v)=H^{\bp}_i(v-)+b^{\bp}_i(x+v)-\left(0\wedge\inf_{u\in[v,v+x]}
(H^{\bp}_i(v-)+b^{\bp}_i(u))\right).$$
Suppose $v$ is such that $H^{\bp}_{m}(v-)=0$ 
(here $H^{\bp}_{m}(0-)=H^{\bp}_{m}(1-)$), call such a number a 
{\em last empty spot}. 
%Notice that such a $v$ might belong to $A_{m}$ despite its name. 
By (ii), Lemma \ref{H}, the set of last empty spots is not empty since it contains
$t_m$. On the other hand, by (i) in the same lemma, the support of $H^{\bp}_m$
is the  closure of $A_m$ which  has measure $1$,  hence it is $\T$. By (iii), we
conclude by letting $v=0,v'\uparrow 1$ that
$H^{\bp}_m(x)=H^{\bp}_m(0)+b^{\bp}_m(x)$ for $0\leq x<1$, so for $x=t_m-$, 
$H^{\bp}_m(0)=-b^{\bp}_m(t_m-)=-\inf b^{\bp}_m$ necessarily since $H^{\bp}_m$
is non-negative. This implies that the last empty spots are those $v$'s such that 
$b^{\bp}_m(v-)=\inf b^{\bp}_m$. We choose one of them by letting
$$V=\inf\{x\in[0,1]:b^{\bp}_m(x-)=\inf_{u\in[0,1]}b^{\bp}_m(u)\},$$ the
first location when the infimum of $b^{\bp}_m$ is reached. We have proved 

\begin{lmm}\label{V}
For any $0\leq x<1,1\leq i\leq m$,
$$H^{\bp}_i(x+V)=b^{\bp}_i(x+V)-\inf_{u\in[V,V+x]}
b^{\bp}_i(u).$$
\end{lmm}

Recall that we are interested in  $\Lambda^{\bp}(i)$, 
the ranked sequence of the lengths of the interval components of $A_{i}$,
where $A_{i}$ can be viewed as the painted portion of the circle after
$i$ drops of paint have fallen, or the set of occupied spots after the $i$-th caravan has arrived.
Lemma \ref{H}(i) enables us to identify $A_{i}$ as the interior of support of 
the function $H^{\bp}_{i}$, and since the Lebesgue measure 
of the interval components of the interior of the 
support of $H^{\bp}_i$ is not affected by a cyclic shift, we record the following 
simple identification

\begin{lmm}\label{Lamb} For every $i=1,\ldots,m$, $\Lambda^{\bp}(i)$
coincides with the ranked lengths of the intervals of constancy of the function
$$x \longmapsto \inf_{u\in[V,V+x]}
b^{\bp}_i(u)\,,\qquad x\in [0,1].$$
\end{lmm}

\section{Convergence of bridges}\label{sec:conv}

We now consider a rescaled randomized version of the bridges introduced above. 
Let $B^{(\eps)}=\eps^{-1+1/\alpha}b^{\bp}_{m}$, where $b^{\bp}_{m}$ is obtained
as above with data $m=T_{\eps},p_{i}=\eps\ell^*_{i},s_{i}=U_{i}$, 
and these quantities are introduced in the Introduction. So for $0\leq x \leq 1$
$$B^{(\eps)}(x)=-\eps^{-1+1/\alpha}x+\sum_{i=1}^{T_{\eps}}\eps^{1/\alpha}
{\ell_i^*}
\ind_{\{x\geq U_i\}}=\eps^{1/\alpha}\sum_{i=1}^{T_{\eps}}
{\ell_i^*}
(\ind_{\{x\geq U_i\}}-x) ,$$
because $\ell^*_1+\ldots+\ell^*_{T_{\eps}}=1/\eps$.
Recall that ${\cal B}^{(2)}$ denotes the standard Brownian bridge, and 
${\cal B}^{(\alpha)}$ the standard stable loop 
with index $\alpha$ as defined in (\ref{pontstable}).

\begin{lmm}
\label{convBeps}
As $\eps\downarrow0$,  the bridge $B^{(\eps)}$  converges weakly
 on the space $\D$  of c\`adl\`ag paths endowed with Skorokhod's topology, to a
bridge with exchangeable increments $B=(B(x),0\leq x\leq 1)$. More precisely:

{\rm (i)} If $\alpha=2$ then $B$ is distributed as $\sqrt{\mu_2/\mu_1}\, {\cal B}^{(2)}$.

{\rm (ii)} If $\alpha\in(1,2)$, then
$B$ is distributed as 
$$\left({\Gamma(2-\alpha)c\over (\alpha-1)\mu_{1}}\right)
^{1\over \alpha}\, {\cal B}^{(\alpha)}\,.$$
\end{lmm}

The proof of Lemma \ref{convBeps}(ii) will use the following well-known representation:
$$\left({\Gamma(2-\alpha)c\over (\alpha-1)\mu_{1}}\right)
^{1\over \alpha}\, 
{\cal B}^{(\alpha)}(x)=
\sum_{i=1}^{\infty}\Delta_i\left(\ind_{\{x\geq U_i\}}-x\right)\, ,
\qquad 0\leq x\leq 1,$$
where  $(U_i,i\geq  1)$  is  a  sequence  of  i.i.d.\  uniform$(0,1)$  r.v.'s,
$(\Delta_i,i\geq 1)$ is the ranked sequence  of the atoms of a Poisson measure
on $(0,\infty)$ with intensity  
$\alpha c\mu_{1}^{-1}x^{-1-\alpha}\d x$, and these
two sequences are independent. 
More precisely, the series in the right-hand side does not converge absolutely,
but is taken in the sense
$$\sum_{i=1}^{\infty}\Delta_i\left(\ind_{\{x\geq U_i\}}-x\right)
=\lim_{n\to\infty}\sum_{i=1}^{n}
\Delta_i\left(\ind_{\{x\geq U_i\}}-x\right),$$
where the limit is uniform in the variable $x$, a.s. This representation
follows immediately from the celebrated L\'evy-It\^{o} decomposition,
specified for the stable process $\sigma^{(\alpha)}$, as the process
of the jumps of  the latter is a Poisson point process on $\R_{+}$
with intensity ${\alpha(\alpha -1)\over \Gamma(2-\alpha)}x^{-1-\alpha}\d x$.
See also Kallenberg \cite{kal73}.

\bigskip

\proof Following Kallenberg \cite{kal73}, 
we represent the jump sizes of the bridge $B^{(\eps)}$
by the random point measure 
$$\psi_{\eps}=\sum_{i=1}^{T_{\eps}}
(\eps^{1/\alpha}\ell_i^*)^2\delta_{\eps^{1/\alpha}\ell_i^*}.$$
By Theorem 2.3 in \cite{kal73}, we have to show:

\begin{equation}\label{i}
\mbox{if $\alpha=2$, then } \psi_{\eps}\to (\mu_2/\mu_1)\delta_0, 
\end{equation}
and
\begin{equation}\label{ii}
\mbox{if $\alpha<2$, then } \psi_{\eps}\to \psi:=
\sum_{i=1}^{\infty}\Delta_i^2\delta_{\Delta_i},
\end{equation}

\noindent where the convergence is in law with respect to the weak
topology on measures on $[0,\infty)$, and in (\ref{ii}), 
$(\Delta_i,i\geq 1)$ is the ranked sequence  of the atoms of a Poisson measure
on $(0,\infty)$ with intensity  
$\alpha c\mu_{1}^{-1}x^{-1-\alpha}\d x$.

Case (i) is easier to treat.  Indeed, notice that the total mass of
$\psi_{\eps}$ is 
$$\psi_{\eps}(\R_+)=\eps\sum_{i=1}^{T_{\eps}}(\ell_i^*)^2
=\eps T_{\eps}\frac{(\ell_{T_{\eps}}^*)^2+\sum_{i=1}^{T_{\eps}-1}\ell_i^2}
{T_{\eps}}.$$
Since $\ell_i^*\leq \ell_i$, the law of large numbers gives 
$\psi_{\eps}(\R_+)\to \mu_2/\mu_1$. 

Now let
$$m^*_{\eps}:=\sqrt{\eps}\max_{1\leq i\leq T_{\eps}}\ell_i^*
\qquad\mbox{and}\qquad M_n:=
\sqrt{\eps}\max_{1\leq i\leq n}\ell_i,$$
so to prove (\ref{i}), it suffices to show that $m^*_{\eps}\to 0$ in 
probability. Notice that $m^*_{\eps}\leq M_{T_\eps}$. 

Let $\eta>0$ and $K>\mu_1^{-1}$. Then 
\begin{eqnarray*}
\P(m^*_{\eps}>\eta)&=&\P(m^*_{\eps}>\eta,T_{\eps}\leq
K\eps^{-1})+\P(m^*_{\eps}>\eta,T_{\eps}>K\eps^{-1})\\
&\leq& \P(M_{\lfloor K\eps^{-1}\rfloor}>\eta)+\P(T_{\eps}>K\eps^{-1}). 
\end{eqnarray*}
The second term converges to $0$ since $\eps\mu_1T_{\eps}\to1$ a.s. 
For the first term, notice that 
$$\P(M_{\lfloor K\eps^{-1}\rfloor}\leq \eta)=
(1-\P(\ell>\eta/\sqrt{\eps}))^{\lfloor K\eps^{-1}\rfloor}.$$
Taking logarithms and checking that $\eps^{-1}\P(\ell^{2}>\eta^{2}/\eps)\to 0$ 
as $\eps\downarrow0$ (which holds since $\E[\ell^2]<\infty$), we finally
obtain that $\P(M_{\lfloor K\eps^{-1}\rfloor}\leq \eta)\to 1$.  
This completes the proof of (\ref{i}).
\medskip

Now we turn our attention to (\ref{ii}). It suffices to show that for every function
$f:[0,\infty)\to [0,\infty)$, say of class ${\cal C}^1$ with bounded derivative
\begin{equation}\label{laplace}
\lim_{\eps\to 0}
\E\left(\exp(-\langle \psi_{\eps},f\rangle)\right) =
\E\left(\exp(-\langle \psi,f\rangle)\right);
\end{equation}
see for instance Section II.3 in Le Gall \cite{legall99}.
In this direction, recall from the classical formula 
for Poisson random measures that
$$\E\left(\exp(-\langle \psi,f\rangle)\right)
={\alpha c\over \mu_{1}} \int_{0}^\infty
(1-\exp (-y^2f(y)))y^{-1-\alpha}dy.
$$

To start with, we observe from the renewal theorem that 
$\eps^{1/\alpha}\ell^*_{T_{\eps}}$ converges to $0$ in probability as
$\eps\to 0$, so in (\ref{laplace}), we may replace $\psi_{\eps}$ by
$$\psi'_{\eps}=\sum_{i=1}^{T_{\eps}-1}
(\eps^{1/\alpha}\ell_i)^2\delta_{\eps^{1/\alpha}\ell_i}.$$

Next, for every $a\geq 0$, we consider the random measure
$$\psi_{\eps,a}=\sum_{i=1}^{a/{\eps}}
(\eps^{1/\alpha}\ell_i)^2\delta_{\eps^{1/\alpha}\ell_i}.$$
Again, by the (elementary) renewal theorem, $\eps T_{\eps}\to \mu_{1}^{-1}$ in probability,
so for every $\eta>0$, the event
\begin{equation}\label{encadrem}
\langle \psi_{\eps,\mu_{1}^{-1}-\eta},f\rangle\,\leq\,
\langle \psi'_{\eps},f\rangle\,\leq \, 
\langle \psi_{\eps,\mu_{1}^{-1}+\eta},f\rangle
\end{equation}\label{encadre}
has a probability which tends to $1$ as $\eps\to 0$.

Now
$$\E\left(\exp(-\langle \psi_{\eps,a},f\rangle)\right)=\E\left(\exp(
-f(\eps^{1/\alpha}\ell)(\eps^{1/\alpha}\ell)^2)\right)^{a/\eps}.$$
Taking logarithms, we have to estimate
\begin{eqnarray*}
& &{a\over \eps}\E\left(1-\exp(-f(\eps^{1/\alpha}\ell)(\eps^{1/\alpha}\ell)^2)
\right)\\
&=&{a\over \eps}\int_{0}^\infty\eps^{2/\alpha}(2xf(\eps^{1/\alpha}x)+
\eps^{1/\alpha}x^2f'(\eps^{1/\alpha}x))
\exp(-(\eps^{1/\alpha}x)^2f(\eps^{1/\alpha}x))\P(\ell>x)dx\\
&=&{a\over \eps}\int_{0}^\infty(2yf(y)+y^2f'(y))
\exp(-y^2f(y))\P(\ell>y/\eps^{1/\alpha})dx.
\end{eqnarray*}
By (\ref{tail}) and dominated convergence, we see that the preceding quantity
converges as $\eps\to 0$ towards
$$ac\int_{0}^\infty(2yf(y)+ y^2f'(y))
\exp(-y^2f(y))y^{-\alpha}dx
=\alpha  ac\int_{0}^\infty(1-\exp(-y^2f(y)))y^{-1-\alpha}dx.$$
Taking $a=\mu_{1}^{-1}\pm\eta$, using (\ref{encadrem}) and letting $\eta$ tend to $0$, we see that (\ref{laplace}) holds, 
which completes the proof of the statement.
 \cq

\section{Convergence of $\bX^{(\eps)}$}\label{convXeps}

In this section, we deduce Theorem \ref{convaddcoal} from Lemmas
\ref{Lamb},\ref{convBeps}. Recall the definition of the bridge $b^{\bp}_i$ in Section \ref{sec:bridge}.
For $i\leq T_{\eps}$, let $B^{(\eps)}_i$ be the bridge 
$\eps^{-1+1/\alpha}b^{\bp}_i$ with data $p_j=\eps \ell^*_{j},s_j=U_j$, so
$B^{(\eps)}_{T_{\eps}}=B^{(\eps)}$. Let also $V_{\eps}$ be the left-most location of
the infimum of $B^{(\eps)}$, and 
$${\bf V}B^{(\eps)}(x)=B^{(\eps)}(x+V_{\eps})-\inf 
B^{(\eps)},\qquad 0\leq x\leq 1$$ 
the Vervaat transform of $B^{(\eps)}$. 
By Lemma \ref{Lamb},  $\bX^{(\eps)}(t)=
\Lambda^{\bp}(T_{\eps}-\lfloor t\eps^{-1/\alpha}\rfloor)$ coincides with the
ranked sequence of lengths of constancy intervals of the infimum process of
$$B^{(\eps)}_{T_{\eps}-\lfloor t\eps^{-1/\alpha}\rfloor}
(x+V_{\eps})-\inf B^{(\eps)},\qquad
0\leq x\leq 1,$$ where the constant $-\inf B^{(\eps)}$ has no effect and is
added for future considerations. 

\begin{lmm}\label{drift}
For every $t\geq 0$, the difference 
$$B^{(\eps)}(x)-B^{(\eps)}_{T_{\eps}-\lfloor
t\eps^{-1/\alpha}\rfloor}(x)
=\eps^{1/\alpha}\sum_{j=0}^{\lfloor t\eps^{-1/\alpha}-1 \rfloor}
\ell^*_{T_{\eps}-i}\ind_{\{x\geq U_{T_{\eps}-i}\}} \qquad 0\leq x\leq 1$$
converges in probability for the uniform norm to the pure drift 
$x\mapsto t\mu_1 x$ as $\eps\downarrow0$. 
\end{lmm}

\proof
Recall from the renewal theorem that $\eps^{1/\alpha}\ell_{T_{\eps}}\to 0$ 
in probability as $\eps\downarrow 0$.  Therefore, we
might start the sum appearing in the statement from $j=1$. Now, the
sequences $(\ell_1,\ldots,\ell_{T_{\eps}-1})$ and
$(\ell_{T_{\eps}-1},\ldots,\ell_1)$ have the same distribution. Up to doing
the substitution, Lemma \ref{drift} for fixed $s$ is therefore a simple application of the
strong law of large numbers.  The conclusion is obtained by standard
monotonicity arguments. \cq

As a consequence of Lemmas \ref{convBeps}, \ref{drift}, and the fact that
$s\mapsto t\mu_1s$ is continuous, the process 
$$B^{(\eps)}_{T_{\eps}-\lfloor t\eps^{-1/\alpha}\rfloor}(x+V_{\eps})-
\inf B^{(\eps)}={\bf V}B^{(\eps)}(x)-\left(B^{(\eps)}(x+V_{\eps})-B^{(\eps)}_{
T_{\eps}-\lfloor t\eps^{-1/\alpha}\rfloor}(x+V_{\eps})\right)$$ 
converges in the Skorokhod space to
$${\cal E}^{(t\mu_1)}=({\cal E}(x)-t\mu_1x, 0\leq x\leq 1),$$ where 
$${\cal E}(x)=B(x+V)-\inf B,\qquad 0\leq x\leq 1$$ is the Vervaat transform of the limiting bridge
$B$ which appears in Lemma \ref{convBeps}, $V$ being  the a.s. unique location of its infimum. 
% Notice that with the
% notation in the Introduction and Sect.\ \ref{sec:conv},
% \begin{equation}\label{E}
% E\build=_{}^{d}\sigma(\Delta)E^{\sigma(\Delta)^{-1}\Delta}. 
% \end{equation}
Now letting $\un{\cal E}^{(t)}$ be the infimum process of ${\cal E}^{(t)}$ and 
$\bF(t)$ be the decreasing sequence of lengths of constancy intervals of 
$\un{\cal E}^{(t)}$, we have 

\begin{prp}\label{convintervalles} 
The process $(\bX^{(\eps)}(t),t\geq 0)$ converges to $(\bF(\mu_1t),t\geq 0)$ in
the sense of weak convergence of finite-dimensional marginals. 
\end{prp}

\proof
The technical point is that Skorokhod convergence of 
$B^{(\eps)}_{T_{\eps}-\lfloor t\eps^{-1/\alpha}\rfloor}(x+V_{\eps})-
\inf B^{(\eps)}$ to ${\cal E}^{(t\mu_1)}$, though it does imply convergence of
respective infimum processes, does not {\em a priori} imply that of the
ranked sequence of lengths of constancy intervals of these processes. However, 
this convergence does hold
%\marginpar{\tiny citer des sources et d\'etailler un peu}
because for every $t\geq 0$, if $(a,b)$ is such a constancy
interval, then ${\cal E}^{(t\mu_1)}(x)>{\cal E}^{(t\mu_1)}(a)$ for
$x\in(a,b)$, a.s. 
%This last fact is obtained
%by standard considerations on L\'evy processes which admit densities. 
See e.g. Lemmas 4 and 7 in \cite{berteac00}.\cq 

This proposition
% according to  the results of \cite{jpda97ebac,berteac00} recalled in the Introduction, 
proves Theorem \ref{convaddcoal}. Indeed, recall from Lemma \ref{convBeps}
that $B=c_{\alpha}{\cal B_{\alpha}}$, where $c_{2}=\sqrt{\mu_{2}/\mu_{1}}$
and for $1<\alpha < 2$
$$c_{\alpha}=\left({\Gamma(2-\alpha)c\over (\alpha-1)\mu_{1}}\right)
^{1\over \alpha}.$$
Then plainly, $\bF(e^{-t})=\bC^{(\alpha)}(t+\log c_{\alpha})$, and hence the limiting process
$\bX(e^{-t})$ is distributed as $\bF(\mu_{1}e^{-t})=\bC^{(\alpha)}(t+\log c_{\alpha}-\log \mu_{1})$.

\section{Related results for a discrete problem}\label{sec:discrete}

In situations involving parking problems, it may be more natural to consider
discrete parking lots, i.e.\ $\Z/n\Z$ instead of the unit circle, and caravans
with integer sizes, e.g.\ as in Knuth's original parking problem.  
Each caravan chooses a random spot, uniform on $\Z/n\Z$,
and tries to park at that spot. Studying the frequencies of blocks of cars 
fits with our general framework by taking
$\ell$ with integer values, $\eps=1/n$ and $s_i=\lfloor n
U_i\rfloor/n$. Rename by $T_n$ the former quantity $T_{\eps}$ (the number of
caravans). Let 
\begin{eqnarray*}
\wt{B}^{(n)}(x)&=&n^{1/\alpha}\sum_{i=1}^{T_n}\left(\frac{\ell^*_i}{n}\ind_{\{x\geq 
\lfloor nU_i\rfloor/n\}}-x\right)\qquad 0\leq x\leq 1\\
B^{(n)}(x)&=&n^{1/\alpha}
\sum_{i=1}^{T_n}\left(\frac{\ell^*_i}{n}\ind_{\{x\geq U_i\}}-x
\right)\qquad 0\leq x\leq 1,
\end{eqnarray*}
so $B^{(n)}$ would be $B^{(1/n)}$ in the  notation above. 
The analog of Lemma \ref{drift} is still true when replacing 
$B^{(n)}$ by $\wt{B}^{(n)}$, without essential change in the
proof. Thus to obtain the very same conclusions as in the preceding sections,
it suffices to check a result similar to Lemma \ref{convBeps}. Namely, we
must prove that $\wt{B}^{(n)}\to B$ in the Skorokhod space as $n\to\infty$. 
Now it is easy to check that a.s., $|\wt{B}^{(n)}(x)-B^{(n)}(\lceil
nx-\rceil/n)|\leq n^{1/\alpha}/n$ for every $n\geq 1, x\in[0,1]$, because 
no $U_i$ is rational a.s. 
Therefore, it suffices to check
that $B^{(n)}(\lceil n\cdot+\rceil/n)$ converges to $B$ in distribution for
the  Skorokhod  topology  on  $\D$.  Up to  using  Skorokhod's  representation
theorem, this is done by taking $f_n=B^{(n)}$ and $\kappa_n(x)=\lceil
nx+\rceil/n$ in the next lemma. 

\begin{lmm}\label{sko}
Let $(f_n,n\geq 1)$ be a sequence of functions converging in $\D$ to
$f$. For $n\in\N$ let also $\kappa_n$ be a right-continuous 
non-decreasing function 
(not necessarily bijective) from $[0,1]$ to $[0,1]$, such
that the sequence $(\kappa_n)$ converges to the identity function uniformly on $[0,1]$. Then $f_n\circ\kappa_n\to f$ in $\D$. 
\end{lmm}

\proof
First consider the case $f_n=f$ for every $n$.  Fix $\eps>0$.  Let
$\kappa_n^{-1}$ be the right-continuous inverse of $\kappa_n$ defined by 
$$\kappa_n^{-1}(x)=\inf\{y\in[0,1]:\kappa_n(y)>x\}.$$ 
It is easy to prove that $\kappa_n(\kappa_n^{-1}(x)-)\leq
x\leq\kappa_n(\kappa_n^{-1}(x))$ for every $x$. 
Since $f$ is c\`adl\`ag, one may find $0=x_0<x_1<\ldots<x_k=1$ such that the
oscillation $\omega(f,[x_i,x_{i+1}))<\eps$ for $0\leq i\leq k-1$, where 
$$\omega(f,A)=\sup_{x,y\in A}|f(x)-f(y)|.$$ 
Since $\kappa_n$ approaches the identity, for $n$ large we may assume 
$\kappa_n(\kappa_n^{-1}(x_i))<\kappa_n(\kappa_n^{-1}(x_{i+1})-)$ for 
$0\leq i\leq k-1$. Define a
time-change $\lambda_n$ (i.e. an increasing bijection between $[0,1]$ and 
$[0,1]$) by interpolating linearly between the points 
$(0,0),(\kappa_n^{-1}(x_i),x_i),1\leq i\leq k-1,(1,1)$. 

Now let
$x\in[0,1]$. Suppose $\kappa_n^{-1}(x_i)\leq x<\kappa_n^{-1}(x_{i+1})$ for some 
$0\leq i\leq k-1$, and notice that $x_i\leq \kappa_n(\kappa_n^{-1}(x_i))\leq
\kappa_n(x)<\kappa_n(\kappa_n^{-1}(x_{i+1})-)\leq x_{i+1}$. Therefore, 
$\kappa_{n}(x)$ belongs to  $[x_{i},x_{i+1})$ as well as $\lambda_{n}(x)$ by definition of 
$\lambda_{n}$, and 
$$|f(\kappa_n(x))-f(\lambda_n(x))|\leq \omega(f,[x_i,x_{i+1}))\leq \eps. $$
Else, one must have $x<\kappa_n^{-1}(0)$ or $x\geq \kappa_n^{-1}(1)$, and the
result is similar. Finally, doing the same reasoning for $\eps=\eps_{n}$ converging to $0$ slowly 
enough gives the existence of some time-changes $\lambda_{n}$ converging to the identity uniformly
such that $\sup_{x\in[0,1]}|f(\kappa_{n}(x))-f(\lambda_{n}(x))|\leq 2\eps_{n}$, hence giving 
convergence of $f\circ\kappa_{n}$ to $f$ in the Skorokhod space. 

In the general case, for every $n\geq 0$ let $\lambda_n$ be a time-change 
such that 
$\lambda_n$ converges to the identity as $n\to\infty$ and 
$f_n\circ\lambda_n$ converges to $f$
uniformly. Take $\kappa'_n=\lambda^{-1}_n\circ \kappa_n$. Then
$f_n\circ\kappa_n-f\circ\kappa'_n\to 0$ uniformly, so it suffices to show that
$f\circ\kappa'_n\to f$ in $\D$, which is done by the former discussion. \cq

In particular, we recover and extend a certain number of results from 
\cite{chasslou99}.

\section{Complements}\label{sec:compl}

In this section, we would like to provide some information 
on  the eternal additive coalescents $\bC^{(\alpha)}$ for $1<\alpha <2$,
which appear in Theorem \ref{convaddcoal}.

\subsection{Mixture of extremes}\label{subsec:extreme}
To start with, we should like to specify the representation of $\bC^{(\alpha)}$ as 
a mixture of so-called {\em extreme}  eternal additive coalescents
 (\cite{jpda97ebac}, \cite{berteac00}).
In this direction, let us first consider a 
sequence $\btheta=(\theta_{0},\theta_{1},\theta_{2},\ldots)$ of non-negative 
numbers satisfying $\sum_{i\geq 0}\theta_i^2=1$ and 
\begin{equation}
\label{condtheta}
\mbox{either }\theta_{0}>0 \qquad\mbox{or}\qquad \sum_{i\geq 0}\theta_{i}=\infty.
\end{equation}
Following Kallenberg \cite{kal73} we associate to $\btheta$ 
a {\em bridge with exchangeable increments}
\begin{equation}
\label{bridgeproc}
B^{\sth}(x)=\theta_{0}\beta(x)+\sum_{i\geq 1}
\theta_{i}(\ind_{\{x\geq U_{i}\}}-x) \qquad 0\leq x\leq 1
\end{equation}
where $(U_{i},i\geq 1)$ denotes a sequence of iid uniform variables and $\beta$ is an independent standard
 Brownian bridge. We write $\bC^{\sth}$ for the eternal additive coalescent
 associated to the bridge $B=B^{\sth}$ as explained in the Introduction and call such 
 $\bC^{\sth}$ extreme.

 According to \cite[Theorem 15]{jpda97ebac},
 every eternal version of the additive coalescent $\bC$ can be obtained as a
mixing of shifted versions of extreme eternal additive coalescents $\bC^{\sth}$,
i.e.\ $\bC$ can be expressed in the form $(\bC^{\sth^*}(t-t^*),t\in\R)$ 
with $\btheta^*,t^*$ random. 
Equivalently, $\bC$ can be viewed as the eternal additive coalescent
constructed in the Introduction from the bridge with exchangeable increments
$B=e^{t^*}B^{\sth^*}$. As observed by Aldous and Pitman 
\cite{jpda97ebac}, the mixing variables $\btheta^*,t^*$ can be recovered
from the initial behavior of $\bC$:
$$e^{t^*}\theta^*_{i}=\lim_{t\to -\infty} e^{-t}\bC_{i}(t)
\quad \hbox{and}\quad e^{2t^*}=\lim_{t\to -\infty} e^{-2t}\sum_{i=1}^{\infty}
\bC_{i}^2(t)\,.$$

In the case of the standard stable loop ${\cal B}^{(\alpha)}$ with $1<\alpha<2$,
recall from the L\'evy-It\^{o} decomposition that 
$\theta_{0}^*=0$ and 
$(e^{t^*}\theta^*_{1}, e^{t^*}\theta^*_{2}, \ldots)=
(\Delta_{1}, \Delta_{2}, \ldots)$
is the ranked sequence of the atoms of a Poisson random measure
on $(0,\infty)$ with intensity ${\alpha(\alpha-1)\over 
\Gamma(2-\alpha)}x^{-1-\alpha}dx$.
In particular, 
\begin{equation}\label{equt}
e^{2t^*}=\sum_{i=1}^{\infty}\Delta_{i}^2
\end{equation}
has the law of a (positive) stable variable with index $\alpha/2$
and
\begin{equation}\label{equth}
\theta^*_{i}=\Delta_{i}/e^{t^*}\,, \qquad i=1,2,\ldots
\end{equation}
is such that  the sequence of squares
$((\theta_{1}^*)^2, (\theta_{2}^*)^2,\ldots)$
is distributed according to the  Poisson-Dirichlet law ${\rm PD}(\alpha/2,0)$;
see Pitman and Yor \cite{py97}.

We also stress that every coalescent $\bC^{\sth}$ can be obtained as a limit of appropriate 
caravan parking problems, which are quite natural given the results of 
\cite{jpda97ebac,berteac00}. Precisely, suppose that a sequence of probabilities 
$\bp^n=(p^n_{1},\ldots,p^n_{m_{n}})$ satisfying 
$p^n_{1}\geq\ldots\geq p^n_{m_{n}}>0$ is given, and satisfies
\begin{equation}\label{regime}
\max_{1\leq i\leq m_{n}}p^n_{i}\build\to_{n\to\infty}^{}0 \qquad \mbox{and}
\qquad
\sigma(\bp^{n})^{-1}p^n_{i}\build\to_{n\to\infty}^{} \theta_{i}\,\qquad 
i\geq1
\end{equation}
for a sequence $\btheta$ as described above, and where 
$\sigma(\bp)=\sqrt{\sum_{i=1}^m p_{i}^2}$ when $\bp=(p_{1},\ldots,p_m)$. For 
every $n$, let $\tau_{n}$ be a uniform permutation on $\{1,2,\ldots,m_{n}\}$.  
Consider the parking problem where the caravans which try to park successively have 
magnitudes $p^n_{\tau_{n}(1)},p^n_{\tau_{n}(2)},\ldots$. Let
$U_{1},U_{2},\ldots$ be independent uniform$(0,1)$ random
variables independent of 
$\tau_{n}$, so we may consider the bridge with exchangeable increments
$$B^{(n)}(s)=\sigma(\bp^{n})^{-1}\left(-x+\sum_{i=1}^{m_{n}}
p^n_{\tau_{n}(i)}\ind_{\{s\geq U_{i}\}}\right)\,, \qquad 0\leq s\leq 1.$$
Kallenberg's theorem shows that under the asymptotic assumptions on $\bp^n$, $B^{(n)}$ 
converges in distribution to the bridge $B^{\sth}$ defined above. 

Now for $t\geq 0$, let 
$I^n_{t}=\inf\{i\geq 1:\sum_{j=i+1}^{m_{n}}p^n_{\tau_{n}(j)}\leq t\}$. The 
following analogue of Lemma \ref{drift} holds.

\begin{lmm}
For every $t\geq 0$, the process 
$$\sigma(\bp^n)^{-1}\sum_{i=I^n_{t}+1}^{m_{n}}p^n_{\tau_{n}(i)}
\ind_{\{s\geq U_{i}\}}\, ,\qquad 0\leq s\leq 1$$
converges in probability for the uniform norm to the pure drift $s\mapsto ts$ as 
$n\to\infty$.
\end{lmm}

\proof
The key to this lemma is to show that 
\begin{equation}\label{maxzero}
\max_{i\geq I^n_{t}}\sigma(\bp^n)^{-1}p^n_{\tau_{n}(i)}\to 0
\end{equation} 
in probability as $n\to\infty$. The result is then obtained via
the so-called ``weak law of large numbers for sampling without replacement'': if 
$x^n_{i},1\leq i\leq n$ is a sequence with sum $t$ satisfying 
$\max_{1\leq i\leq n} x^n_{i}\to0$ as $n\to\infty$, and if $\tau_{n}$ is a uniform
permutation on $\{1,\ldots,n\}$, then for every rational $r\in[0,1]$,
$\sum_{i=1}^nx^n_{\tau_{n}(i)}\ind_{\{r\geq U_{i}\}}\to tr$ in probability (in fact in 
$L^2$). The result in probability remains true if $x^n_{i},1\leq i\leq n$ is random 
with sum $t$, and $\max_{1\leq i\leq n} x^n_{i}\to 0$ in probability. 
One concludes that the 
process $(\sum_{i=1}^nx^n_{\tau_{n}(i)}\ind_{\{s\geq U_{i}\}},0\leq s\leq 1)$ converges 
in probability to $(ts,0\leq s\leq 1)$ for the uniform norm by a monotonicity argument. 
The lemma is then proved by letting $x_{1}=\sigma(\bp^n)^{-1}p^n_{\tau_{n}(I^n_{t}+1)},
x_{2}=\sigma(\bp^n)^{-1}p^n_{\tau_{n}(I_{t}+2)},
\ldots,x_{m_{n}-I^n_{t}}=\sigma(\bp^n)^{-1}p^n_{\tau_{n}(m_{n})},
x_{m_{n}-I^n_{t}+1}=t-\sum_{i=I^n_{t}+1}^{m_{n}}
\sigma(\bp^n)^{-1}p^n_{\tau_{n}(I^n_{t}+1)}$ (note that this last term is 
$\leq \sigma(\bp^n)^{-1}p^n_{\tau_{n}(I^n_{t})}$, which goes to $0$).

So let us show (\ref{maxzero}). To this end, let $0<\rho<1$, then
$X^{\rho}_{n}:=\sum_{i=\lfloor \rho m_{n}\rfloor}^{m_{n}}\sigma(\bp^n)^{-1}
p^n_{\tau_{n}(i)}\to \infty$ in probability, since 
$E[X^{\rho}_{n}]\sim\sigma(\bp^n)^{-1}(1-\rho)$ goes to infinity (notice $\sigma(\bp)\leq p_{1}$)
while $E[(X^{\rho}_{n})^2]\sim E[X^{\rho}_{n}]^2$, as a simple computation shows. 
Therefore, $I^n_{t}\sim m_{n}$ in probability. Consequently, for any $K\in\N$, the quantity
$P(\tau_{n}^{-1}(1)< I^n_{t},\ldots,\tau_{n}^{-1}(K)< I^n_{t})$ goes to $1$, so 
$\min_{i\geq I^n_{t}}\tau_{n}(i)\to\infty$ in probability. But then, for any $\eps>0$, 
if $K$ is such that $\theta_{K}<\eps/2$, then $\sigma(\bp^n)^{-1}p^n_{K}\leq\eps$
for $n$ large. Up to taking $n$ even larger, 
with probability close to $1$, $\tau_{n}(i)\geq K$ for $i\geq I^n_{t}$ and therefore 
$\max_{i\geq I^n_{t}}\sigma(\bp^n)^{-1}p^n_{\tau_{n}(i)}\leq \eps$, hence 
(\ref{maxzero}). \cq

One deduces, as around the proof of Proposition \ref{convintervalles}, the following claim. 
Let $\bX^{(n)}(t)=\Lambda^{\bp^n\circ\tau_{n}}(I^n_{t})$ be as above with data
$m=m_n,p^n_{\tau_{n}(i)},1\leq i\leq m_{n},s_i=U_i$. Then 

\begin{prp}\label{thetacoal}
As $n\to\infty$, under the asymptotic regime (\ref{regime}), the process
$(\bX^{(n)}(t),t\geq 0)$ converges in the sense of weak convergence of finite-dimensional 
marginals to the time-reversed eternal additive coalescent 
$(\bC^{\sth}(-\log(t)),t\geq 0)$. 
\end{prp}

\subsection{On the marginal distributions}
It would also be interesting to determine the marginal laws of the fragmentation
$\bF^{(\alpha)}(t):=\bC^{(\alpha)}(-\log t)$. 
The task seems quite difficult if started from the description
of $\bF^{(\alpha)}(t)$ in terms of lengths of constancy intervals of Vervaat transform of bridges,
because excursion theory seems powerless here, unlike in
\cite{mier01}. In particular, the fact that the fragmentation is based 
on stable loops and not stable bridges impedes the application 
of results of Miermont \cite{mier01}  on additive coalescents 
based on bridges of certain L\' evy processes.

Another way to start the exploration is to use the
representation of  fragmentation processes $\bF^{\sth}(t)
:=\bC^{\sth}(-\log t)$ described  in the preceding section with the
help of Inhomogeneous Continuum Random Trees (ICRT) discussed in
\cite{jpda97ebac}. In particular,  it is easy to obtain the  first moment of a
size-biased pick\footnote{Recall that a size-biased pick $X_{\dagger}$ from
a (random) positive sequence $(X_i,i\geq 1)$ with sum $0<S<\infty$ a.s.\ 
is a random variable of the form $X_{i^*}$, where 
$P(i^*=i|X_j,j\geq 1)=X_i/S$.} 
$\bF_{\dagger}(t)$ from the sequence $\bF(t)$ 
for any fixed $t$, as follows. 

Let us recall the basic facts on the ICRT$(\btheta)$ construction of 
$\bF^{\sth}$. The ICRT can be viewed via a 
{\em stick-breaking construction} as the metric completion of the 
positive real line $\R_+$ endowed with a non standard metric. Precisely, 
suppose we are given the following independent random elements:
\begin{itemize}
\item A Poisson process $\{(U_i,V_i),i\geq 1\}$ on the octant 
$\mathbb{O}=\{(x,y):0<y<x\}\subset \R_+^2$, with intensity 
$\theta_0\d x\d y\ind_{\mathbb{O}}$, so in particular 
$\{U_i,i\geq 1\}$ is a Poisson process with intensity 
$\theta_0 x\d x\ind_{x\geq 0}$,
\item A sequence of independent Poisson processes 
$\{\xi_{i,j},j\geq 1\}, i=1,2,\ldots$ with respective intensities 
$\theta_i\d x\ind_{x\geq 0},i=1,2,\ldots$. 
\end{itemize}
We distinguish the points $(V_i,i\geq 1),(\xi_{i,1},i\geq 1)$ as 
{\em joinpoints}, while $(U_i,i\geq 1),(\xi_{i,j},i\geq 1,j\geq 2)$ are 
called {\em cutpoints}.  If $\eta$ is a cutpoint, let $\eta^*$ be its
associated joinpoint, i.e. $U_i^*=V_i,\xi_{i,j}^*=\xi_{i,1}$. 
By the assumption on $\btheta$, it is a.s.\ possible to arrange the 
cutpoints by increasing order $0<\eta_1<\eta_2<\ldots$. We then construct 
a family ${\cal R}(k),k\geq 1$ of ``reduced trees'' as follows. Cut 
the set $(0,\infty)$ into ``branches'' $(\eta_i,\eta_{i+1}]$, where 
by  convention $\eta_0=0$.  Let  ${\cal R}(1)$  be  the segment  $(0,\eta_1]$,
endowed
with the usual distance $d_1(x,y)=|x-y|$. 
Then given ${\cal R}(k),d_k$, we obtain ${\cal R}(k+1)$
by adding the branch $(\eta_{k},\eta_{k+1}]$ somewhere on ${\cal R}(k)$, 
and we plant the left-end $\eta_k$ on the joinpoint $\eta^*_k$ (since 
a.s.\ $\eta^*<\eta$, the point $\eta_k^*$ is indeed an element of 
${\cal R}(k)$). Precisely, ${\cal R}(k+1)=(0,\eta_{k+1}]$ and 
$d_{k+1}(x,y)=d_k(x,y)$ if $x,y\in{\cal R}(k)$, 
$d_{k+1}(x,y)=|x-y|$ if $x,y\in(\eta_k,\eta_{k+1}]$, and
$d_{k+1}(x,y)=x-\eta_k+d_k(y,\eta_k^*)$ if 
$x\in(\eta_k,\eta_{k+1}],y\in{\cal R}(k)$. 
As the distances $d_k$ are compatible by definition, this defines a
random metric space $(0,\infty),d$ such that the restriction of $d$ to 
${\cal R}(k)$ is $d_k$, we call its metric completion $\TT^{\sth}$ the 
ICRT$(\btheta)$, its elements are called {\em vertices}. 
The point $\varnothing=\lim_{n\to\infty}1/n$ is distinguished
and called the {\em root}. 

One can see that $\TT^{\sth}$ is an $\R$-tree, i.e.\ a complete metric space
such that for any $x,y\in\TT^{\sth}$ there is a unique simple path $[[x,y]]$ 
from $x$ to $y$, which is isometric to the segment $[0,d(x,y)]$, i.e. is a 
geodesic. Moreover, it can be endowed with a natural measure $\mu^{\sth}$ 
which is the weak limit as $n\to\infty$ of the empirical measures 
$n^{-1}\sum_{i=1}^n\delta_{\eta_i}$.  This measure is non-atomic and 
supported on {\em leaves}, i.e.\ vertices $x\in\TT^{\sth}$ such that $x\notin
[[\varnothing,y]]\setminus\{y\}$ for any $y\in\TT^{\sth}$. Non-leaf vertices
form a set called the {\em skeleton}. A second natural measure is the 
Lebesgue measure $\lambda$ on $\TT^{\sth}$, i.e.\ the unique measure such 
that $\lambda([[x,y]])=d(x,y)$ for any $x,y$, and this measure is supported
on the skeleton.

Now for each $t$ consider a Poisson measure on $\TT^{\sth}$ 
with atoms $\{x_i^t,i\geq 1\}$, with intensity $t\lambda(\d x)$, so the 
different processes are coupled in the natural way as $t$ varies, i.e.\ 
$\{x_i^t,i\geq 1\}$ increases with $t$. These points disconnect the 
tree into a forest of disjoint connected tree components, order them
as $\FF^{\sth}_i(t),i\geq 1$ by decreasing order of $\mu^{\sth}$-mass. Then
the process $((\mu^{\sth}(\FF^{\sth}_i(t)),i\geq 1),t\geq 0)$ of these 
$\mu^{\sth}$-masses has same law as $\bF^{\sth}$. A size-biased pick from 
this sequence of masses is then obtained as the $\mu^{\sth}$-mass of the 
tree component at time $t$ that contains an independent
$\mu^{\sth}$-sample, conditionally on $(\TT^{\sth},\mu^{\sth})$. Therefore, 
if $\bF^{\sth}_{\dagger}(t)$ denotes such a size-biased pick, 
$E[\bF^{\sth}_{\dagger}(t)]$ is the probability that two independent 
$\mu^{\sth}$-samples $X_1,X_2$ belong to the same tree component of the cut 
tree, i.e.\ that no atom of the Poisson measure at time $t$ falls in the path 
$[[X_1,X_2]]$, and hence it equals $E[e^{-td(X_1,X_2)}]$. 

It turns out \cite{jpda97ebac} that $d(X_1,X_2)$ has same law as the 
length $\eta_1$ of the first branch (i.e.\ the length of ${\cal R}(1)$).
It is easy to see (see also \cite{jpmc97b}) that this branch's length 
has law 
$$P(\eta_1>r)=e^{-\theta_0^2r2/2}\prod_{i=1}^{\infty}
(1+\theta_ir)e^{-\theta_ir}.$$

In our setting, recall that the random sequence $\theta^*$ is related to that 
of the atoms $(\Delta_{i})$ of a Poisson measure on $(0,\infty)$ with intensity
${\alpha(\alpha-1)\over \Gamma(2-\alpha)}x^{-1-\alpha}dx$
by (\ref{equt}) and (\ref{equth}).
Observe that $\bF^{(\alpha)}(t)=\bF^{\theta^*}(t/e^{t^*})$, and since
we must take a Poisson process with intensity $t/e^{t^*}$ on the skeleton of the 
ICRT($\theta^*$), terms $e^{t^*}$ cancel out and 
$E[\bF_{\dagger}^{(\alpha)}(t)]=E[e^{-t\eta}]$ where
\begin{eqnarray*}
P(\eta\geq r)&=&E\left[\prod_{i=1}^{\infty}(1+r\Delta_i)
e^{-r\Delta_i}\right]\\
&=&\exp\left(-\int_0^{\infty}\frac{\alpha(\alpha-1)\d x}{\Gamma(2-\alpha)
x^{1+\alpha}}
(1-\exp(-rx+\log(1+rx)))\right)\\
&=&\exp(-(\alpha-1)r^{\alpha}),
\end{eqnarray*}
which is a Weibull distribution. This gives (at least in principle) the first
moment 
$$E[\bF_{\dagger}^{(\alpha)}(t)]=
\int_0^{\infty}\alpha(\alpha-1)r^{\alpha-1}\exp\left(-t r-(\alpha-1)
r^{\alpha}\right)\d r.$$

In principle, this method could be used for the computation of moments of higher order, 
where the length $\eta_{1}$ of the first branch would simply be replaced by the 
total length $\eta_{k}$ of ${\cal R}(k)$. Unfortunately, the distribution of 
$\eta_{k}$ is complicated for $k=3$, and seems intractable for higher $k$'s 
(\cite{jpda99tree}).

\def\polhk#1{\setbox0=\hbox{#1}{\ooalign{\hidewidth
  \lower1.5ex\hbox{`}\hidewidth\crcr\unhbox0}}}

\end{document}

%% file: park.pstex_t
\begin{picture}(0,0)%
\includegraphics{park.ps}%
\end{picture}%
\setlength{\unitlength}{3947sp}%
\begingroup\makeatletter\ifx\SetFigFont\undefined%
\gdef\SetFigFont#1#2#3#4#5{%
  \reset@font\fontsize{#1}{#2pt}%
  \fontfamily{#3}\fontseries{#4}\fontshape{#5}%
  \selectfont}%
\fi\endgroup%
\begin{picture}(7544,3514)(3129,-6344)
\put(4351,-5836){\makebox(0,0)[lb]{\smash{{\SetFigFont{12}{14.4}{\rmdefault}{\mddefault}{\updefault}{\color[rgb]{0,0,0}$A_i$}%
}}}}
\put(5176,-6286){\makebox(0,0)[lb]{\smash{{\SetFigFont{12}{14.4}{\rmdefault}{\mddefault}{\updefault}{\color[rgb]{0,0,0}$s_{i+1}$}%
}}}}
\put(8026,-6286){\makebox(0,0)[lb]{\smash{{\SetFigFont{12}{14.4}{\rmdefault}{\mddefault}{\updefault}{\color[rgb]{0,0,0}$t_{i+1}$}%
}}}}
\put(6226,-5836){\makebox(0,0)[lb]{\smash{{\SetFigFont{12}{14.4}{\rmdefault}{\mddefault}{\updefault}{\color[rgb]{0,0,0}$A_i$}%
}}}}
\put(7351,-5836){\makebox(0,0)[lb]{\smash{{\SetFigFont{12}{14.4}{\rmdefault}{\mddefault}{\updefault}{\color[rgb]{0,0,0}$A_i$}%
}}}}
\put(8851,-5836){\makebox(0,0)[lb]{\smash{{\SetFigFont{12}{14.4}{\rmdefault}{\mddefault}{\updefault}{\color[rgb]{0,0,0}$A_i$}%
}}}}
\put(5551,-3961){\makebox(0,0)[lb]{\smash{{\SetFigFont{12}{14.4}{\rmdefault}{\mddefault}{\updefault}{\color[rgb]{0,0,0}$i+1$-th caravan}%
}}}}
\put(5851,-2986){\makebox(0,0)[lb]{\smash{{\SetFigFont{12}{14.4}{\rmdefault}{\mddefault}{\updefault}{\color[rgb]{0,0,0}$p_{i+1}$}%
}}}}
\end{picture}%

%% file: park2.pstex_t
\begin{picture}(0,0)%
\includegraphics{park2.ps}%
\end{picture}%
\setlength{\unitlength}{3947sp}%
\begingroup\makeatletter\ifx\SetFigFont\undefined%
\gdef\SetFigFont#1#2#3#4#5{%
  \reset@font\fontsize{#1}{#2pt}%
  \fontfamily{#3}\fontseries{#4}\fontshape{#5}%
  \selectfont}%
\fi\endgroup%
\begin{picture}(7524,2445)(3139,-6673)
\put(9151,-5311){\makebox(0,0)[lb]{\smash{{\SetFigFont{12}{14.4}{\rmdefault}{\mddefault}{\updefault}{\color[rgb]{0,0,0}$H^{\bp}_i$}%
}}}}
\put(6226,-6136){\makebox(0,0)[lb]{\smash{{\SetFigFont{12}{14.4}{\rmdefault}{\mddefault}{\updefault}{\color[rgb]{0,0,0}$A_i$}%
}}}}
\put(7351,-6136){\makebox(0,0)[lb]{\smash{{\SetFigFont{12}{14.4}{\rmdefault}{\mddefault}{\updefault}{\color[rgb]{0,0,0}$A_i$}%
}}}}
\put(8851,-6136){\makebox(0,0)[lb]{\smash{{\SetFigFont{12}{14.4}{\rmdefault}{\mddefault}{\updefault}{\color[rgb]{0,0,0}$A_i$}%
}}}}
\put(4351,-6136){\makebox(0,0)[lb]{\smash{{\SetFigFont{12}{14.4}{\rmdefault}{\mddefault}{\updefault}{\color[rgb]{0,0,0}$A_i$}%
}}}}
\put(4876,-5086){\makebox(0,0)[lb]{\smash{{\SetFigFont{12}{14.4}{\rmdefault}{\mddefault}{\updefault}{\color[rgb]{0,0,0}$p_{i+1}$}%
}}}}
\put(5176,-6136){\makebox(0,0)[lb]{\smash{{\SetFigFont{12}{14.4}{\rmdefault}{\mddefault}{\updefault}{\color[rgb]{0,0,0}$s_{i+1}$}%
}}}}
\put(8026,-6136){\makebox(0,0)[lb]{\smash{{\SetFigFont{12}{14.4}{\rmdefault}{\mddefault}{\updefault}{\color[rgb]{0,0,0}$t_{i+1}$}%
}}}}
\put(6301,-4936){\makebox(0,0)[lb]{\smash{{\SetFigFont{12}{14.4}{\rmdefault}{\mddefault}{\updefault}{\color[rgb]{0,0,0}$h^{\bp}_{i+1}$}%
}}}}
\end{picture}%